\documentclass[a4paper,12pt]{amsart}
\usepackage{amssymb}
\usepackage{ifthen}
\usepackage{graphicx}
\usepackage{float}
\usepackage{caption}
\usepackage{subcaption}
\usepackage{cite}
\usepackage{amsfonts}
\usepackage{amscd}
\usepackage{amsxtra}
\usepackage{color}

\newtheorem{thm}{Theorem}
\newtheorem{cor}{Corollary}
\newtheorem{lem}{Lemma}
\newtheorem{rem}{Remark}

\newtheorem{conj}{Conjecture}
\newtheorem{prob}{Problem}

\theoremstyle{definition}
\newtheorem{defn}{Definition}[section]
\newtheorem{example}{Example}

\newenvironment{pf}[1][]{%
 \vskip 1mm
 \noindent
 \ifthenelse{\equal{#1}{}}%
  {{\slshape Proof. }}%
  {{\slshape #1.} }%
 }%
{\qed\bigskip}

\newcounter{alphabet}
\newcounter{tmp}
\newenvironment{Thm}[1][]{\refstepcounter{alphabet}%
\bigskip%
\noindent%
{\bf Theorem \Alph{alphabet}}%
\ifthenelse{\equal{#1}{}}{}{ (#1)}%
{\bf .} \itshape}{\vskip 8pt}

\makeatletter
\newcommand{\Ref}[1]{\@ifundefined{r@#1}{}{\setcounter{tmp}{\ref{#1}}\Alph{tmp}}}
\makeatother

\newenvironment{Lem}[1][]{\refstepcounter{alphabet}%
\bigskip%
\noindent%
{\bf Lemma \Alph{alphabet}}%
{\bf .} \itshape}{\vskip 8pt}

\newcommand{\IC}{{\mathbb C}}
\newcommand{\ID}{{\mathbb D}}

\def\be{\begin{equation}}
\def\ee{\end{equation}}

\newcommand{\bee}{\begin{enumerate}}
\newcommand{\eee}{\end{enumerate}}

\newcommand{\blem}{\begin{lem}}
\newcommand{\elem}{\end{lem}}
\newcommand{\bthm}{\begin{thm}}
\newcommand{\ethm}{\end{thm}}
\newcommand{\bcor}{\begin{cor}}
\newcommand{\ecor}{\end{cor}}
\newcommand{\beg}{\begin{example}}
\newcommand{\eeg}{\end{example}}
\newcommand{\begs}{\begin{examples}}
\newcommand{\eegs}{\end{examples}}
\newcommand{\bdefe}{\begin{defn}}
\newcommand{\edefe}{\end{defn}}
\newcommand{\bprob}{\begin{prob}}
\newcommand{\eprob}{\end{prob}}
\newcommand{\bques}{\begin{ques}}
\newcommand{\eques}{\end{ques}}
\newcommand{\bei}{\begin{itemize}}
\newcommand{\eei}{\end{itemize}}
\newcommand{\bcon}{\begin{conj}}
\newcommand{\econ}{\end{conj}}
\newcommand{\bcons}{\begin{conjs}}
\newcommand{\econs}{\end{conjs}}
\newcommand{\bprop}{\begin{propo}}
\newcommand{\eprop}{\end{propo}}
\newcommand{\br}{\begin{rem}}
\newcommand{\er}{\end{rem}}
\newcommand{\brs}{\begin{rems}}
\newcommand{\ers}{\end{rems}}
\newcommand{\bo}{\begin{obser}}
\newcommand{\eo}{\end{obser}}
\newcommand{\bos}{\begin{obsers}}
\newcommand{\eos}{\end{obsers}}
\newcommand{\bpf}{\begin{pf}}
\newcommand{\epf}{\end{pf}}
\newcommand{\ba}{\begin{array}}
\newcommand{\ea}{\end{array}}
\newcommand{\beq}{\begin{eqnarray}}
\newcommand{\beqq}{\begin{eqnarray*}}
\newcommand{\eeq}{\end{eqnarray}}
\newcommand{\eeqq}{\end{eqnarray*}}

\newcommand{\ds}{\displaystyle}

\newcounter{minutes}
\divide\time by 60
\newcounter{hours}
\multiply\time by 60 \addtocounter{minutes}{-\time}

\begin{document}
\bibliographystyle{amsplain}
\title[Sharp Bohr type inequality]{Sharp Bohr type inequality}

\thanks{
File:~\jobname .tex,
          printed: \number\day-\number\month-\number\year,
          \thehours.\ifnum\theminutes<10{0}\fi\theminutes}


\author{Amir Ismagilov}
\address{A. Ismagilov, Kazan Federal University, Kremlevskaya 18, 420 008 Kazan, Russia}
\email{amir.ismagilov@list.ru}

\author{Ilgiz R. Kayumov}
\address{I. R. Kayumov, Kazan Federal University, Kremlevskaya 18, 420 008 Kazan, Russia}
\email{ikayumov@kpfu.ru}

\author{Saminathan Ponnusamy}
\address{S. Ponnusamy, Department of Mathematics, Indian Institute of
Technology Madras, Chennai-600 036, India. }

\email{samy@iitm.ac.in}



\begin{abstract}
This article is devoted  to the sharp improvement of the classical Bohr inequality for bounded analytic functions defined on the unit disk.
We also prove two other sharp versions of the Bohr inequality by replacing the constant term by the absolute of the function and
the square of the absolute of the function, respectively.
\end{abstract}


\maketitle
\pagestyle{myheadings}
\markboth{Amir Ismagilov,  Ilgiz R Kayumov and Saminathan Ponnusamy}{Bohr type inequalities}

\section{Introduction and Main Results}\label{KayPon8-sec1}
Let $\ID :=\{z\in\IC:\, |z|<1\}$ be the open unit disk in the complex plane $\IC$.
Bohr's theorem (after subsequent improvements due to M.~Riesz, I.~Schur and F.~Wiener) states that if
$\sum_{k=0}^\infty a_k z^k$ is a bounded analytic function in $\ID$  such that  $|f(z)| \leq 1$ in $\ID$, then \cite{Bohr-14} \be\label{ClassBohr-1}
B_f(r):=\sum_{k=0}^\infty |a_k|r^k\leq 1 ~\mbox{ for any $r\leq 1/3$},
\ee
and the value $1/3$ is sharp. There are many proofs of this inequality (cf. \cite{Sidon-27-15} and \cite{Tomic-62-16}).
 Also, we would to recall that Bombieri \cite{Bom-62} proved the following
which gives an upper bound on the growth of $B_f(r)$:
$$
\sup B_f(r)
=\frac{3-\sqrt{8(1-r^2)}}{r}, \quad \mbox{$1/3\leq  r\leq 1/\sqrt{2}$}.
$$
Bombieri and Bourgain \cite{BomBor-04} investigated asymptotical behaviour of Bohr's sums as $r \to 1$. In fact,
they constructed $a_n$, and by a delicate analysis of exponential sums, proved that when $r\rightarrow 1$,
$$ B_f(r) \geq (1-r^2)^{-1/2} - \left( c \log \frac{1}{1-r} \right)^{3/2+\epsilon},
$$
where $c=c(\epsilon)$ depends on $\epsilon$.

This result has created enormous interest on Bohr's inequality in various settings.
See for example, \cite{BoasKhavin-97-4,AlKayPON-19, BhowDas-18, BenDahKha,BomBor-04,DF,DFOOS,DGM,DeGarM04,DMS,DP,KayPon1,KayPon2,KayPon3,KayPon3c,PauSin--09,Pop}
and the recent survey on this topic by
Abu-Muhanna et al. \cite{AAPon1},  \cite{IKKayPon1},  \cite[Chapter 8]{GarMasRoss-2018} and the references
therein. Some of these articles use various methods from complex analysis, functional analysis, number theory, and probability,
and furthermore they provide new theory and  applications of Bohr's results on his work on Dirichlet series.
 For example, several multidimensional generalizations of this result are obtained in
\cite{Aizen-00-1,aiz,AAD,Aiz07,BoasKhavin-97-4,DjaRaman-2000}. Moreover, in \cite{KP-AASFM2-19}, the authors
answer the open problem about the powered Bohr radius posed by Djakov and Ramanujan  in 2000.

For the case $a_0=0$,  Tomi\'c \cite{Tomic-62-16} showed that \eqref{ClassBohr-1} holds for $0\leq r\leq 1/2$.
Later Ricci \cite{Ricci-1956} improved it by showing that \eqref{ClassBohr-1} holds for $0\leq r\leq 3/5$,
and the largest value of $r$ for which \eqref{ClassBohr-1} holds in this case would be in the interval $(3/5,1/\sqrt{2}]$. Later in 1962, Bombieri \cite{Bom-62}  obtained
the sharp upper estimate for $\sum_{k=0}^\infty |a_k|r^k$ in the case $r\in (1/3,1/\sqrt{2}]$.
See \cite{KayPon1,KayPon2,KP-AASFM2-19} for new proofs of it in a general form.

If $f$ has a higher order zero at the origin, then we have only a partial answer about the range for $r$
(see \cite[Remark 2]{PaulPopeSingh-02-10} and also \cite{PonWirths-20} for some related investigation
on this problem).

It is worth pointing out that the notion of Bohr's radius, initially defined for analytic functions from the unit disk
$\ID$ into $\ID$, was generalized by authors to include mappings from $\ID$ to some other domains $\Omega$ in $\ID$ (\cite{Abu,Abu2,Aiz07}).

Unless otherwise stated, throughout this paper  $S_r(f)$ denotes the area of the
image of the subdisk $|z|<r$ under the mapping $f$ and when there is no confusion, we let for brevity $S_r$ for $S_r(f)$.
More recently,  Kayumov and Ponnusamy  \cite{KayPon3b} improved  Bohr's inequality \eqref{ClassBohr-1} in various forms. For example,
the following inequality is obtained in \cite{KayPon3} (see also \cite{KayPon3b}).

\begin{Thm}\label{KayPon8-Additional4}
Suppose that $f(z) = \sum_{k=0}^\infty a_k z^k$ is analytic in $\ID$ and $|f(z)| \leq 1$ in $\ID$. Then
\begin{equation}\label{Eq_Th3}
|a_0|+\sum_{k=1}^\infty |a_k|r^k+\frac{16}{9}\left (\frac{S_r}{\pi}\right )  \leq 1 ~\mbox{ for  }~ r \leq \frac{1}{3},
\end{equation}
and the numbers $1/3$ and $16/9$ cannot be improved.

\end{Thm}
%
%

Also, it is worth recalling that if the first term $|a_0|$ in \eqref{Eq_Th3} is replaced by $|a_0|^2$, then  $1/3$ and $16/9$ could be replaced
by $1/2$ and $9/8$, respectively (cf. \cite{KayPon3}). For a harmonic analog of Bohr's inequality \eqref{ClassBohr-1} and
that of Theorem \Ref{KayPon8-Additional4}, we refer to \cite{EPR-2017} and \cite{KayPon3}, respectively. More recently,  the
authors in \cite{IKKayPon1} improved Theorem \Ref{KayPon8-Additional4} by replacing the quantity $S_r/\pi$  in the
inequality \eqref{Eq_Th3} by $S_r/(\pi-S_r)$. Moreover,  in \cite{KayPon3b} (see also \cite[Corollary 3]{IKKayPon1}),
the following generalizations of Bohr's results were obtained. The idea used in these two articles \cite{IKKayPon1}, \cite{KayPon1} were
used for further investigation on Bohr's inequality by Liu et al.  \cite{LSX2018}.

\begin{Thm}\label{KayPon8-cor1}
Suppose that $f(z) = \sum_{k=0}^\infty a_k z^k$ is analytic in $\ID$ and $|f(z)| <1$ in $\ID$. Then
$$|f(z)|+\sum_{k=1}^\infty |a_k|r^k \leq 1 ~\mbox{ for  }~ |z|=r \leq \sqrt{5}-2 \approx 0.236068
$$
and no larger radius than $\sqrt{5}-2 $ will do. Moreover,
$$|f(z)|^2+\sum_{k=1}^\infty |a_k|r^k \leq 1 ~\mbox{ for  }~ |z|=r  \leq 1/3
$$
and no larger radius than $1/3$ will do.
\end{Thm}

Our main goal of this article is to derive sharp version of Theorems \Ref{KayPon8-Additional4} and \Ref{KayPon8-cor1}, and hence
certain results from the survey article \cite{IKKayPon1}. It is important to point out that for individual functions the Bohr radius is always
greater than $1/3$.

We now state our main results and their proofs will be presented in the next section.
Moreover, in the interesting results which are presented below, there is an extremal function
such that $1/3$ cannot be increased.  The price for such important fact is essentially due to
nonlinearity of the obtained functionals.

\begin{thm}\label{VINITI-1}
Suppose that $f(z) = \sum_{k=0}^\infty a_k z^k$ is analytic in $\ID$ and $|f(z)| \leq 1$ in $\ID$.
Then
\begin{equation}\label{EqVINITI-1}
M(r):=\sum_{k=0}^\infty |a_k|r^k+\frac{16}{9}\left (\frac{S_r}{\pi}\right )+\lambda\left (\frac{S_r}{\pi}\right )^2  \leq 1 ~\mbox{ for  }~ r \leq \frac{1}{3}
\end{equation}
where
 $$
 \lambda= \frac{4 (486 - 261 a - 324 a^2 + 2 a^3 + 30 a^4 + 3 a^5)}{ 81 (1 + a)^3 (3 - 5 a)}=18.6095 \ldots
 $$ and  $a\approx 0.567284$, is the unique positive root of the equation $\psi (t)=0$ in the interval $(0,1)$, where
$$ \psi (t)= -405 + 473 t + 402 t^2 + 38 t^3 + 3 t^4 + t^5.
$$
 The equality is achieved for the function
$$f(z) = \frac{a-z}{1-az}.
$$

\end{thm}

\br It is evident that $\lambda$ is an algebraic number. Moreover, it can be shown that
$\lambda$ is indeed the positive root of the algebraic equation
$$
285212672 + 6268596224 x + 37178714880 x^2 + 87178893840 x^3 +
$$
$$
+ 97745285925 x^4 - 5509980288 x^5=0.
$$
\er

\br
From the proof of Theorem \ref{VINITI-1}, the following observation is clear: for any function $F:[0,\infty ) \to [0,\infty)$ such that $F(t)>0$ for $t>0$,
there exist analytic functions $f:\, \mathbb{D} \to \mathbb{D}$ for which the inequality
$$\sum_{k=0}^\infty |a_k|r^k+\frac{16}{9}\left (\frac{S_r}{\pi}\right )+\lambda\left (\frac{S_r}{\pi}\right )^2 +F(S_r) \leq 1
~\mbox{ for  }~ r \leq \frac{1}{3}
$$
is false, where $\lambda$ is as in Theorem \ref{VINITI-1}.
\er This result follows from the fact that there is a concrete function for which the equality in (\ref{EqVINITI-1}) holds so that one can add no  strictly
 positive terms.

%
%

\begin{thm}\label{VINITI-2}
Suppose that $f(z) = \sum_{k=0}^\infty a_k z^k$ is analytic in $\ID$ and $|f(z)| \leq 1$ in $\ID$.
Then
\begin{equation}\label{EqVINITI-2}
M(z,r):=|f(z)|^2+\sum_{k=1}^\infty |a_k|r^k+\frac{16}{9}\left (\frac{S_r}{\pi}\right )+\lambda\left (\frac{S_r}{\pi}\right )^2  \leq 1 ~\mbox{ for  }~ |z| = r \leq \frac{1}{3},
\end{equation} where
$$\lambda=\frac{-81 + 1044 a + 54 a^2 - 116 a^3 - 5 a^4}{162 (a+1)^2 (2 a-1)}=16.4618 \ldots
$$ and $a\approx 0.537869$ is the unique positive root of the equation
$$
-513 + 910 t + 80 t^2 + 2 t^3 + t^4=0
$$
in the interval $(0,1)$.

The equality is achieved for the function
$$f(z) = \frac{a-z}{1-az}.
$$

\end{thm}
\br One can check  that actually
$\lambda$ is the positive root of the algebraic equation
$$
575930368 + 4437874624 x + 11353360788 x^2 + 10868034060 x^3 -  703096443 x^4=0.
$$
\er

 One can replace $|f(z)|^2$ by $|f(z)|$ in Theorem \ref{VINITI-2} but this will decrease the Bohr radius. Namely, the following theorem is valid.

\begin{thm}\label{Additional}
Suppose that $f(z) = \sum_{k=0}^\infty a_k z^k$ is analytic in $\ID$ and $|f(z)| \leq 1$ in $\ID$.
Then
\begin{equation}\label{EqVINITI-3a}
|f(z)|
+ \sum_{k=1}^\infty |a_k| r^k + p\left (\frac{S_r}{\pi}\right ) \leq 1 ~\mbox{ for  }~  |z| =r \leq \sqrt{5}-2 \approx 0.236068,
\end{equation}
where the constants  $r_0= \sqrt{5} -2$ and $p=2 (\sqrt{5}-1)$ are sharp.
\end{thm}

\section{Proofs of Theorem \ref{VINITI-1}, \ref{VINITI-2} and \ref{Additional}}

For the proof of our results, we need the following lemmas.

\begin{Lem}\label{KP2-lem2}
{\rm \cite[Lemma 2]{KayPon1}}
 Let $|b_0|<1$ and $0 < r \leq 1$. If $g(z)=\sum_{k=0}^{\infty} b_kz^k$ is analytic and satisfies the inequality $|g(z)| \leq 1$ in $\ID$, then
the following sharp inequality holds:
\begin{equation}\label{KP2-eq3a1}
\sum_{k=1}^\infty |b_k|^2r^{pk} \leq r^{p}\frac{(1-|b_0|^2)^2}{1-|b_0|^2r^{p}}.
\end{equation}
\end{Lem}

\begin{Lem}\label{KP2-lem2a1}
{\rm \cite[Lemma 1]{KayPon3b}}
Suppose that $g(z) = \sum_{k=0}^\infty b_k z^k$ is analytic in $\ID$, $|g(z)|<1$ in $\ID$ and  $S_r(g)$ denotes the area of the
image of the subdisk $|z|<r$ under the mapping $g$. Then
the following sharp inequality holds:
\begin{equation}\label{KP2-eq3}
\frac{S_r(g)}{\pi}:= \sum_{k=1}^\infty k |b_k|^2r^{2k} \leq  r^2\frac{(1-|b_0|^2)^2}{(1-|b_0|^2r^2)^2}
~\mbox{ for $0 < r \leq 1/\sqrt{2}$}.
\end{equation}
\end{Lem}

\br Lemma \Ref{KP2-lem2} is not true for $r>1/\sqrt{2}$  unless $g$ is univalent in $\ID$. For instance, consider $g(z)=z^n$ where $n\geq 2$.
\er

\begin{Lem}\label{KP2-lem2b}
Let $p \in \mathbb{N}$, $0\leq m \leq p$ and $f(z)=\sum_{k=0}^{\infty} a_{pk+m}z^{pk+m}$ be analytic in $\ID$
and $|f(z)| < 1$ in $\ID$. Then the following inequalities hold:
\be\label{KayPon8-eq7a}
\sum_{k=1}^{\infty} |a_{pk+m}|r^{pk} \leq
\left \{ \begin{array}{lr} \ds 
r^p\frac{1-|a_m|^2}{1-r^p|a_m|} & \mbox{ for $|a_m| \ge r$} \\[4mm]
\ds 
r^p\frac{\sqrt{1-|a_m|^2}}{\sqrt{1-r^{2p}}} &\mbox{ for $|a_m| < r$}.
\end{array} \right .
\ee
\end{Lem}
\bpf
Proof of this lemma follows from the proof of Theorem 1 in \cite{KayPon2}  (see also \cite{IKKayPon1}, and \cite[Proof of Theorem ~1]{KayPon1}
for the case $m=0$ and \cite{KayPon2}). The proof uses the classical Cauchy-Schwarz inequality and \eqref{KP2-eq3a1}.
However, for the sake of completeness, we supply here some details.  At first, we write
$f(z)=z^m g(z^p)$, where $|g(z)| \leq 1$ in $\ID$ and
$g(z)=\sum_{k=0}^{\infty} b_kz^k$ is analytic in $\ID$  with $b_k=a_{pk+m}$.  Let $|b_0|=a$. Choose any $\rho >1$ such that $\rho r \leq 1$. Then it follows from the classical Cauchy-Schwarz inequality and \eqref{KP2-eq3a1} with $\rho r$ in place of $r$ that
\begin{eqnarray*}\label{Additional1}
\sum_{k=1}^{\infty} |a_{pk+m}|r^{pk} &= &
\sum_{k=1}^{\infty}|b_k| r^{pk} \\
&\leq &  \sqrt{\sum_{k=1}^{\infty}|b_k|^2 \rho^{pk}r^{pk}} \sqrt{\sum_{k=1}^{\infty} \rho^{-pk}r^{pk}} \\
&\leq &
\sqrt{r^p\rho^p\frac{(1-a^2)^2}{1-a^2r^p\rho^p}}\, \sqrt{\frac{\rho^{-p}r^p}{1-\rho^{-p}r^{p}}} ~\mbox{ (by Lemma C)}\\
&=& \frac{r^p(1-a^2)}{\sqrt{1-a^2r^p\rho^p}}\,  \frac{1}{\sqrt{1-\rho^{-p}r^{p}}}.
\end{eqnarray*}
We need to consider the cases $a \ge r^p$ and $a<r^p$, independently. If $a \ge r^p$, then we may set $\rho=1/\sqrt[p]{a}$.
On the other hand, if $a < r^p$, then we set $\rho=1/r$. As a result of these substitutions, we easily obtain that
\begin{equation*}
\sum_{k=1}^{\infty}|b_k| r^{pk}  \leq  \left \{
\begin{array}{lr}
\ds r^p\frac{(1-a^2)}{1-r^p a} & \mbox{ for $a \ge r^p$}\\
\ds r^p \frac{\sqrt{1-a^2}}{\sqrt{1-r^{2p}}} & \mbox{ for $a <r^p$}
\end{array} \right .
\end{equation*}
and the desired inequalities \eqref{KayPon8-eq7a} follow by setting  $b_k=a_{pk+m}$ and $a=|a_m|$.
\epf

In the proof of our main results, we just need Lemma \Ref{KP2-lem2a1} and the case $m=0$ and $p=1$ of Lemma \Ref{KP2-lem2b}.

\vspace{8pt}

\subsection{Proof of Theorem \ref{VINITI-1}}
Consider the function $M(r)$ given  by \eqref{EqVINITI-1}. Since $M(r)$ is an increasing function of $r$,
we have $M(r)\leq M(1/3)$ for $0\leq r\leq 1/3$ and thus,
it suffices to prove the inequality \eqref{EqVINITI-1} for $r=1/3$.
Moreover,  for $m=0$ and $p=1$,  Lemma \Ref{KP2-lem2b} gives
\be\label{KayPon8-eq7}
\sum_{k=1}^\infty |a_k|r^k \leq
\left \{ \begin{array}{lr} \ds A(r):=r\frac{1-|a_0|^2}{1-r|a_0|} & \mbox{ for $|a_0| \ge r$} \\[4mm]
\ds B(r):=r\frac{\sqrt{1-|a_0|^2}}{\sqrt{1-r^2}} &\mbox{ for $|a_0| < r$}.
\end{array} \right .
\ee
At first we consider the case $|a_0| \ge 1/3$. In this case, using \eqref{KayPon8-eq7} (with $r=1/3$) and Lemma \Ref{KP2-lem2a1}, we have
\beqq
M(1/3) &\leq&  |a_0|+A(1/3) +\frac{16}{9\pi} S_{1/3}+\lambda\left (\frac{S_{1/3}}{\pi}\right )^2\\
&\leq & |a_0|+ \frac{1-|a_0|^2}{3-|a_0|}+16 \frac{(1-|a_0|^2)^2}{(9-|a_0|^2)^2}+81 \lambda \frac{(1-|a_0|^2)^4}{(9-|a_0|^2)^4}\\
& =&1- \frac{(1-|a_0|)^3}{(9-|a_0|^2)^4} \Phi(|a_0|),
\eeqq
where
\beqq
\Phi(t)&=& 3078 + 1944 t - 522 t^2 - 432 t^3 + 2 t^4 + 24 t^5 + 2 t^6  \\
   && ~~~~+\lambda(-81  - 243 t  - 162 t^2  + 162t^3  + 243 t^4  + 81 t^5 ).
\eeqq
One can verify that the function $\Phi(t)$ in the interval $[1/3,1]$ has exactly one stationary point $t_0=0.567284\ldots $ which
is the positive root of the equation $ \Phi' (t)=0$.

Let us show that $t_0=a$ and that $\Phi(t_0)=0$. The equation $\Phi'(a)=0$ is fulfilled automatically (in fact, $\lambda$ was chosen in such the way at the beginning). A little algebra gives
$$
\Phi(a)=2\frac{a^2-9}{3-5a} \psi(a).
$$

Consequently, $\Phi(a)=0$. Besides this observation, we have  $\Phi(1/3)>0$ and $\Phi(1)>0$. Thus, $\Phi(t)\geq 0$ in the interval $[1/3,1]$
which proves that $M(r)\leq 1$ for $|a_0|\in [1/3,1]$ and for $r\leq 1/3$

Next we consider the case $|a_0|<1/3$. Again,  using \eqref{KayPon8-eq7} (with $r=1/3$) and Lemma \Ref{KP2-lem2a1}, we deduce that
\beqq
M(1/3) &\leq & |a_0|+B(1/3) +16 \frac{(1-|a_0|^2)^2}{(9-|a_0|^2)^2}+81 \lambda \frac{(1-|a_0|^2)^4}{(9-|a_0|^2)^4}\\
& \leq& |a_0|+\frac{\sqrt{1-|a_0|^2}}{\sqrt{8}}+16 \frac{(1-|a_0|^2)^2}{(9-|a_0|^2)^2}+81 \lambda \frac{(1-|a_0|^2)^4}{(9-|a_0|^2)^4}\\
&=& \Psi (|a_0|).
\eeqq
Routine and straightforward calculations show that the last expression for $\Psi (t)$ is an increasing function of $t$
and so its maximum is achieved at the point $t_0=|a_0|=1/3$, and the maximum value of $\Psi (t_0)$ is seen to be less than $1$
(in fact it is $\leq 0.98$) so that the desired inequality \eqref{EqVINITI-1} follows for $|a_0|\in [0,1/3)$.
This proves that $M(r)\leq 1$ for $|a_0|\in [0,1/3)$ and for $r\leq 1/3$.

To prove that the constant $\lambda$ is sharp, we consider the function $f$ given by
\be\label{EqSharpTh1}
f(z) = \frac{a-z}{1-az} =a - (1-a^2)\sum_{k=1}^\infty a^{k-1} z^{k}, \quad z\in\ID,
\ee
where $a\in (0,1)$, and compute the value of $M(r)$ for this function. Indeed, we may let
$$ M_1(1/3)= \sum_{k=0}^\infty |a_k|3^{-k} + \frac{16}{9\pi} S_{1/3}+\lambda_1\left (\frac{S_{1/3}}{\pi}\right )^2,
$$
where $a_0=a$ and $a_k= -(1-a^2) a^{k-1}$. Straightforward calculations show that
\beqq
M_1(1/3)
&=& a+ \frac{1-a^2}{3-a}+16 \frac{(1-a^2)^2}{(9-a^2)^2}+81 \lambda \frac{(1-a^2)^4}{(9-a^2)^4}+81(\lambda_1-\lambda)\frac{(1-a^2)^4}{(9-a^2)^4}
\eeqq
Choose $a$ as the positive root $t_0$ of the equation $\psi (t)=0$.
As a consequence, we see that
$$ M_1(1/3)= 1+81(\lambda_1-\lambda)\frac{(1-a^2)^4}{(9-a^2)^4}
$$
which is bigger than $1$ in case $\lambda_1> \lambda$. This proves the sharpness assertion and the proof is complete.



\subsection{Proof of Theorem \ref{VINITI-2}}
Let $M(z,r)$ be defined by \eqref{EqVINITI-2}. As $M(z,r)$ is an increasing function of $r$,
it suffices to prove the inequality \eqref{EqVINITI-2} for $r=1/3$. Moreover, by the assumption and the
Schwarz-Pick lemma applied to the function $f$ show that
$$|f(z)|\leq \frac{r+|a_0|}{1+r|a_0|}=:D(r), \quad |z| \leq r.
$$
At first we consider the case $|a_0| \ge 1/3$. In this case, using \eqref{KayPon8-eq7} (with $r= 1/3$), Lemma \Ref{KP2-lem2a1} and the last inequality, we have
for $r\leq 1/3$ that
\beqq
M(z,r) &\leq& M(z,1/3)\\
&\leq&  \left (D(1/3)\right )^2 +A(1/3) +\frac{16}{9\pi} S_{1/3}+\lambda\left (\frac{S_{1/3}}{\pi}\right )^2\\
&\leq &  1- \left [ 1-\left (\frac{1+3|a_0|}{3+|a_0|}\right )^2
- \frac{1-|a_0|^2}{3-|a_0|}-16 \frac{(1-|a_0|^2)^2}{(9-|a_0|^2)^2}-81 \lambda \frac{(1-|a_0|^2)^4}{(9-|a_0|^2)^4}\right ]\\
&= &  1- (1-|a_0|)^2(1+|a_0|)\left [ \frac{|a_0|+15}{(3+|a_0|)^2(3-|a_0|)} -16 \frac{(1+|a_0|)}{(9-|a_0|^2)^2}\right .\\
&& \hspace{5cm} \left .-81 \lambda \frac{(1-|a_0|^2)^2(1+|a_0|)}{(9-|a_0|^2)^4}\right ]\\
&= &  1- (1-|a_0|)^3(1+|a_0|)\left [ \frac{|a_0|+29}{(9-|a_0|^2)^2 } -81 \lambda \frac{(1-|a_0|^2)(1+|a_0|)^2}{(9-|a_0|^2)^4}\right ]\\
& =&1- \frac{(1-|a_0|)^3(1+|a_0|)}{(9-|a_0|^2)^4} \Phi(|a_0|),
\eeqq
where
\beqq
\Phi(t)=2349 + 81 t - 522 t^2 - 18 t^3 + 29 t^4 + t^5 + \lambda(-81 - 162 t +
 162 t^3 + 81 t^4), ~ 1/3\leq t\leq 1.
\eeqq

We see that the function $\Phi(t)$ in the interval $[1/3,1]$ has exactly one stationary point $t_0=0.537869\ldots $ which is the positive
root of the equation $\Phi'(t)=0$.

Let us show that $t_0=a$ and that $\Phi(t_0)=0$. The equation $\Phi'(a)=0$ is fulfilled automatically. A little algebra gives
$$
\Phi(a)=\frac{9-a^2}{2(2a-1)}( -513 + 910 a + 80 a^2 + 2 a^3 + a^4)
$$
so that $\Phi(t_0)=0$.

Besides this observation,  we find that $\Phi(1/3)>0$ and $\Phi(1)>0$. Thus, $\Phi(t)\geq 0$ in the interval $[1/3,1]$
which proves that $M(z,r)\leq 1$ for $|a_0|\in [1/3,1]$ and $r\leq 1/3$. Thus, the desired inequality \eqref{EqVINITI-2} follows for
$|a_0|\in [1/3,1]$.

Next we consider the case $|a_0|<1/3$. Again,  using \eqref{KayPon8-eq7} (with $r= 1/3$) and Lemma \Ref{KP2-lem2a1}, we deduce that
\beqq
M(z,r) &\leq &   \left (D(1/3)\right )^2 +B(1/3) +\frac{16}{9\pi} S_{1/3}+\lambda\left (\frac{S_{1/3}}{\pi}\right )^2\\
& \leq& \left (\frac{1+3|a_0|}{3+|a_0|}\right )^2+\frac{\sqrt{1-|a_0|^2}}{\sqrt{8}}+16 \frac{(1-|a_0|^2)^2}{(9-|a_0|^2)^2}+81 \lambda \frac{(1-|a_0|^2)^4}{(9-|a_0|^2)^4}\\
&=& \Psi (|a_0|).
\eeqq
Routine and straightforward calculations show that the last expression for $\Psi (t)$ is an increasing function of $t$ in the interval $[0,1/3]$,
and so its maximum is achieved at the point $t_0=|a_0|=1/3$, and the maximum value of $\Psi (t_0)$ is seen to be less than $1$
(in fact it is $\leq 0.987$) so that the desired inequality \eqref{EqVINITI-2} follows for $|a_0|\in [0,1/3)$. This proves
the inequality \eqref{EqVINITI-2}.

Finally, to prove that the constant $\lambda$ is sharp, as in the previous theorem, we consider the function $f$ given by \eqref{EqSharpTh1}
and compute the value of $M(z,r)$ for this function. Indeed, we may let
$$ M_1(z,r)= \left (\frac{1+3a}{3+a}\right )^2+\sum_{k=1}^\infty |a_k|3^{-k} + \frac{16}{9\pi} S_{1/3}+\lambda_1\left (\frac{S_{1/3}}{\pi}\right )^2,
$$
where $a_0=a$ and $a_k= -(1-a^2) a^{k-1}$.
Straightforward calculations show that
\beqq
M_1(z,r)
&=&  \left (\frac{1+3a}{3+a}\right )^2+ \frac{1-a^2}{3-a}+16 \frac{(1-a^2)^2}{(9-a^2)^2}+81 \lambda \frac{(1-a^2)^4}{(9-a^2)^4}
+81(\lambda_1-\lambda)\frac{(1-a^2)^4}{(9-a^2)^4}.
\eeqq
Choose $a$ as the positive root $t_0$ of the equation $-513 + 910 t + 80 t^2 + 2 t^3 + t^4=0$.
As a consequence, we see that
$$ M_1(z,r))= 1+81(\lambda_1-\lambda)\frac{(1-a^2)^4}{(9-a^2)^4}
$$
which is bigger than $1$ when $\lambda_1> \lambda$.  The proof of the theorem is complete.

\subsection{Proof of Theorem \ref{Additional}}
Let
$$M(r) := \sup_{|z|=r}|f(z)| + \sum_{k=1}^\infty |a_k| r^k + p\left (\frac{S_r}{\pi}\right ).
$$
In the case $|a_0|=a \ge r$ from (\ref{KayPon8-eq7}) it follows that
$$M(r) \leq M_1(r)=\frac{r+a}{1+ar} + r\frac{1-a^2}{1-ra} + p\frac{(1-a^2)r^2}{(1-a^2r^2)^2}.
$$
Clearly, it suffices to show that
$$M_1(\sqrt{5}-2) = \frac{a+\sqrt{5}-2}{\left(\sqrt{5}-2\right) a+1}
+ \frac{\left(\sqrt{5}-2\right) \left(a^2-1\right)}{\left(\sqrt{5}-2\right) a-1}
+ p\frac{\left(9-4 \sqrt{5}\right) \left(a^2-1\right)^2}{\left(\left(4 \sqrt{5}-9\right) a^2+1\right)^2} \leq 1
$$
which is equivalent to the following inequality
$$\frac{(1-a)^3 (7 (-9 + 4 \sqrt{5}) + 4 (-47 + 21 \sqrt{5}) a + (-161 + 72 \sqrt{5}) a^2)}{\left(\left(4 \sqrt{5}-9\right) a^2+1\right)^2} \leq 0.
$$
The last inequality is clearly valid in the unit interval $[0,1]$ because $-9 + 4 \sqrt{5} \leq 0$,
$-47 + 21 \sqrt{5} \leq 0$ and $-161 + 72 \sqrt{5} \leq 0$.

It means that  the desired inequality \eqref{EqVINITI-3a} is correct in the case $|a_0| \geq r$ and $r \leq \sqrt{5}-2$.

Next we consider the $|a_0|=a \leq r$. In this case,  from (\ref{KayPon8-eq7}) we obtain that
$$ M(r) \leq M_2(r)= \frac{r+a}{1+ar} +r\frac{\sqrt{1-a^2}}{\sqrt{1-r^2}} + p\frac{(1-a^2)r^2}{(1-a^2r^2)^2}.
$$
Again routine computations show that the last expression is an increasing function of $r$ for $a \in [0,r]$ and
hence its maximum is achieved at the point $a=r$. This case was settled in the previous considerations.

Let us finally show that the estimate for $p$ cannot be improved. We set
$$f(z)=\frac{z+a}{1+az}
$$ and choose $z=r=\sqrt{5}-2$. One can check that
$$|f(z)| + \sum_{k=1}^\infty |a_k| r^k + (2 (\sqrt{5}-1)+\varepsilon) \frac{S_r}{\pi} =$$
$$\frac{a+\sqrt{5}-2}{\left(\sqrt{5}-2\right) a+1} + \frac{\left(\sqrt{5}-2\right) \left(a^2-1\right)}{\left(\sqrt{5}-2\right) a-1} + (2 (\sqrt{5}-1)+\varepsilon)\frac{\left(9-4 \sqrt{5}\right) \left(a^2-1\right)^2}{\left(\left(4 \sqrt{5}-9\right) a^2+1\right)^2} = $$

$$1+ \frac{(1-a)^3 \left(a \left(\left(72 \sqrt{5}-161\right) a+84 \sqrt{5}-188\right)+7 \left(4 \sqrt{5}-9\right)\right)}{\left(\left(4 \sqrt{5}-9\right) a^2+1\right)^2}+ \varepsilon\frac{\left(9-4 \sqrt{5}\right) \left(a^2-1\right)^2 }{\left(\left(4 \sqrt{5}-9\right) a^2+1\right)^2}.$$
From here we see that in the case when $a\rightarrow 1$ this expression behaves like $1+C(1-a)^2\varepsilon>1$. This concludes the proof.


\subsection*{Acknowledgements}
We thank the referee for his careful reading of our paper and his proposals that helped to ameliorate it.
The work of A. Ismagilov and I. Kayumov is supported by the Russian Science Foundation under grant 18-11-00115.
The  work of the third author is supported by Mathematical Research Impact Centric Support of Department of 
Science \& Technology, India  (MTR/2017/000367).


\end{document}